\documentclass{article}

\newcounter{formula}
\newcommand{\tag}[1]{\refstepcounter{formula}
  \eqno{(\arabic{formula})}\label{#1}}

\newcommand{\NNN}{{\cal N}}
\newcommand{\MMM}{{\cal M}}
\newcommand{\ZZmn}{\ZZ^{m\times n}}
\newcommand{\Conv}{\mathop{\rm Conv}}
\newcommand{\Cone}{\mathop{\rm Cone}}
\newcommand{\CharCone}{\mathop{\rm CharCone}}

\newcommand{\mod}{\mathop{\rm mod}}

\newcommand{\ZZ}{{\bf Z }}
\newcommand{\RR}{{\bf R }}
\newcommand{\QQ}{{\bf Q }}

\newcommand{\Sj}{\mathop{\sum}\limits_{j=1}^{n}}
\newcommand{\Sum}{\mathop{\sum}\limits}

\newcommand{\set}[1]{\left \{ #1 \right\}}
\newcommand{\antie}[1]{\left \lfloor #1 \right\rfloor}
\newcommand{\bl}{\bigl(}
\newcommand{\br}{\bigr)}

\newtheorem{theorem}{Theorem}

\newtheorem{note}{Note}

\newenvironment{remark}{\begin{note} \rm}{\end{note}}

\begin{document}

\title{On the number of vertices \protect\\ in integer linear programming
problems\protect\\{\normalsize Overview}}

\author{Nikolai Yu. Zolotykh\thanks{University of Nizhni Novgorod,
Russia {\tt http://www.uic.nnov.ru/\~{}zny}.}}

\date{1997, last revision: October 25, 2000}

\maketitle

\begin{abstract}
We give a survey of work on the number of vertices of the convex
hull of integer points defined by the system of linear
inequalities. Also, we present our improvement of some of these.
\end{abstract}


\section*{Introduction}
\addcontentsline{toc}{section}{Introduction}

In the integer linear programming problem, it is required to maximize the
linear function $\Sj c_j x_j$, subject to $\Sj a_{ij} x_j \le b_i$ $(i=1 , \dots ,
n)$, $x_j \in \ZZ$ $(j= 1 , \dots , n)$ where $a_{ij}$, $c_j$, $b_i$  are integers.
We can formulate this in the following matrix form:
$$
\max cx \tag{max}
$$
$$
Ax \le b,       \tag{Axleb}
$$
$$
x \in \ZZ^n,    \tag{xinZ}
$$
where $A=(a_{ij}) \in \ZZmn$, $b=(b_i)\in \ZZ^m$, $c=(c_j) \in \ZZ^n$. The set
of all solutions to (\ref{Axleb}) is a polyhedron $P \subseteq \RR^n$.
Consider the set $M(A, b) =P \cap \ZZ^n$ and its convex hull $P_I$.  It is
known that $P_I$ is a polyhedron.  If the maximum (\ref{max}) over the set
$M(A, b)$ associated with (\ref{Axleb},\,\ref{xinZ}) is finite then it is
attained at some vertex of $P_I$; on the other hand, the facets of $P_I$ are
the strongest regular cuts.  Therefore examining a characterization of
$P_I$ can indicate ways of generating effective algorithms or understand why
constructing such algorithms are impossible~\cite{EKK,Schrijver,book}.

Let $N(A, b)$ be the set of vertices (extreme points) of the polyhedron $P_I$.
Here we give a survey of work on the cardinality $|N(A, b)|$.
Also, we give one improvement (see Theorem~\ref{myimprovement}).

It was apparently Rubin~\cite{Rubin} who first showed that there is no function
$g(n, m)$ such that for any $A \in \ZZmn$, $b \in \ZZ^m$ it holds $|N(A, b)|
\le g(n, m)$. In the papers considered below, upper bounds on $|N(A, b)|$ are
(as a rule) functions of three arguments $n, m, \alpha$ where $\alpha = \max
\set{|a_{ij} (i=1 , \dots , m; j=1 , \dots , n)}$.  Another way of upper bounds
determination is to use the length of the expression for the coefficients in
$Ax \le b$. Both ways are equivalent in the following sense. Let
$$\varphi = \max \set{1 + \Sj \bl 1+ \log (1+|a_{ij}|) \br + \log (1+
|b_i|)},$$ then $\varphi \le (1+n)\bl 1+ \log(1+ \alpha)\br $ and, on the other
hand, $\log (1+ \alpha) \le \varphi$.  Hence, if for some function $g(\cdot,
\cdot, \cdot)$ it holds that $|N(A, b)| \le g(n, m, \varphi)$, then $|N(A, b)|
\le g\Bigl(n, m, (1+n)\bl 1+ \log (1+ \alpha) \br \Bigr)$; on the other hand,
if for some function $f(\cdot, \cdot, \cdot)$ it holds that $|N(A, b)| \le
f\bigl(n, m, \log(1+ \alpha)\bigr)$, then $|N(A, b)| \le f(n, m, \varphi)$.

Together with observation of the number of vertices in $P_I$
associated with an arbitrary system $Ax \le b$, we consider the following
important case of the problem. Let  $a_0, a_1 , \dots , a_n$ be natural numbers,
$a=(a_1 , \dots , a_n)$, $\alpha = \max \set{a_0 , \dots , a_n}$. Denote by $\MMM
(a,a_0)$ the set of integer non-negative solutions to the inequality $\Sj a_j
x_j \le a_0$, let $\NNN (a,a_0)$ be the set of vertices of the convex hull of
$\MMM (a,a_0)$. The convex hull of $\MMM (a,a_0)$ is called the {\em knapsack
polytope}.

We note that Chapter~3 in~\cite{book} is specially devoted to topics
considered here. We consider also other results.

Below $\Conv\set{v_1 , \dots , v_l}$ and $\Cone\set{v_1 , \dots , v_l}$
are respectively the convex hull of vectors $v_1 , \dots , v_l$ in
$\RR^n$ and the cone generated by these; $\antie{\alpha}$ is the maximal
integer not exceeding the number $\alpha$ in $\RR$, $\log\alpha$ denotes
the logarithm of $\alpha$ to the base $2$.

\section{Upper bounds on the number of vertices}

An approach of Shevchenko~\cite{1981} has been widely used for
determination of upper bounds on the number of vertices in integer
linear programming problems. It was developed in~\cite{1984,1992,342,book}.
The following property lies at the basis of it.
We say that a set $X \subset \ZZ_+^n$ has the
{\em separation} property if for any distinct points
$x=(x_1 , \dots , x_n)$ and $y=(y_1 , \dots , y_n)$ in $X$ there is an
index $j$ such that $2y_j < x_j$. It has been proved in~\cite{1981} that if
$\max \set{x_j , x\in X} \le k_j -1$ $(j=1 ,\dots , n-1)$ then
$$
|X| \le \prod _{j=1}^{n-1} \Bigl( 1+ \log (k_j -1) \Bigr).  \tag{Razd}
$$
Shevchenko's approach consists of the following steps:
\begin{enumerate}
 \item a one-to-one mapping of the set
       $N(A, b)$ into the set $X$ with separation property is constructed;
 \item quantities $k_j$ are estimated and a bound
       (\ref{Razd}) is written out.
\end{enumerate}

As shown in~\cite{1981}, the set $N'(A,b)$ of all vertices of the convex
hull of $\set{x \in \ZZ ^n, Ax = b, x \ge 0}$ has the separation
property. This and (\ref{Razd}) lead us to the bound
$$
|N'(A,b)| \le \Bigl( 1+\log n + r \log (\alpha \sqrt{r})
\Bigr) ^{r-1}, \tag{N'}
$$
where $r$ is the rank of $A$ and $\alpha$ is the maximum of the
moduli of coefficients in the system $Ax=b$.

This result can be applied to an arbitrary system $Ax
\le b$. Using the standard substitution $x_j ' = \max \set{x_j , 0}$,
$x_j '' = -\min \set{x_j , 0}$, $y=A_0 - Ax$ we lead the system $Ax \le
b$ to the form
$$ \begin{array}{c} Ax'-Ax''+y=a_0,\\[0.25cm]
x' \in \ZZ^n,\, x'' \in \ZZ^n,\, y \in \ZZ^m, \\[0.25cm] x' \ge 0,\, x'' \ge
0,\, y\ge 0.  \end{array} \tag{priv}
$$
It is easy to verify that different points in $N(A, b)$ pass into
different vertices of the convex hull of the set of solutions to the system
(\ref{priv}). From (\ref{N'}) we get the bound~\cite{1981}:
$$
|N(A,b)| \le \bigl( 1 + \log (n+1) + n \log (\alpha \sqrt{n}) \bigr)^{2n+m-1}.
$$
For fixed $n$ and $m$ it has the form of polynomial in $\log \alpha$.
This bound was subsequently improved in~\cite{Cook,342,ChirkovMA2},
and took the form of some polynomial in $m$ and $\log \alpha$ for any
fixed $n$ (see below).

For the knapsack problem we introduce the variable
$x_{n+1}=a_0 - \Sj a_j x_j$ and use the bound (\ref{Razd}).
Then we get~\cite{1981}
$$
|\NNN(a,a_0) \le \prod _{j=1} ^n \left(1+ \log\left(1+\frac{a_0}{a_j}\right)
\right).
\tag{knapsack1}
$$

The use of this approach for a square system of linear inequalities see
in~\cite{1984}.

One more approach to estimate the number of vertices in integer
linear programming problems was proposed by Hayes and Larman \cite{HayesLarman}
and developed by Morgan \cite{Morgan} and Cook et al. \cite{Cook}. It can be
considered as a variant of Shevchenko's method. Let $P$ be a polyhedron,
defined by the system $Ax \le b$.  We divide $P$ into {\em reflecting sets}
$P_1 , \dots , P_l$ such that for any two points $x, y$ in $\ZZ ^n \cap P_i$
either $2x-y$ or $2y-x$ belongs to $P_i$. Then it is not hard to verify that no
reflecting set in the divizion of $P$ contains more than one point in
$N(A, b)$. Thus, $|N(A, b)| \le l$. For obtaining a bound it is enough to
divide $P$ into the ``small'' number of reflecting sets.

Let us consider the use of this approach~\cite{HayesLarman} for estimating
the number of vertices in the knapsack polytope.
The polytope
$$
P= \set{x\ge 0, \Sj
a_j x_j \le a_0}
$$
is divided into boxes
$$
\set{x=(x_1 , \dots , x_n) \in \RR^n , 2^{i_j -1} \le x_j < 2^{i_j}, j=1 , \dots , n},
\tag{Box}
$$
where $i_j \in \ZZ$, $1 \le i_j \le 1 + \log \left( \frac{a_0}{a_j} + 1\right)$
$(j=1 , \dots , n)$. It is obvious that no box contains more than one
point in $\NNN (a, a_0)$. This again leads us to the bound (\ref{knapsack1}).

As Morgan~\cite{Morgan} note,
some boxes (\ref{Box}) can fail to contain points in $\NNN(a,
a_0)$. The reason is that they are ``too far'' from the facets of $P$.
For the estimate of $|\NNN(a, a_0)|$ it is enough to bound the number of
boxes (\ref{Box}) which intersect with the hyperplain $\Sj \frac{x_j}{c_j}
= 1$, where $c_j = \antie{\frac{a_0}{a_j}}$ $(j=1 , \dots , n)$. It is
not hard to prove~\cite{Morgan} that for $n \ge 2$ the number of such
boxes is not exceed the quantity
$$
n \log(2n) \left( 1+ \log \left(1 + \frac{a_0}{\gamma} \right) \right)^{n-1},
$$
where $\gamma = \min \set{|a_j|, j=1 , \dots , n}$.
Thus, we have
\begin{theorem} \label{tMorgan} \cite{Morgan}
For $n \ge 2$ it holds that
$$
|\NNN(a,a_0)|
\le n \log(2n) \left( 1+ \log \left(1 + \frac{a_0}{\gamma} \right)
\right)^{n-1}.  \tag{knapsack2}
$$
\end{theorem}

Cook et al.\,\cite{Cook} generalized this approach for the estimate of the
number of vertices $N(A,b)$ if $A \in \QQ ^{m\times n}$ and $b \in \QQ ^{m}$
are arbitrary. If the length of the expression for the each inequality of $Ax
\le b$ is not more than $\varphi$, then $|N(A,b)|\le 2m^n (6n^2 \varphi )
^{n-1}$.  For systems with integer coefficients we have $\varphi \le
(n+1)\bigl(1+ \log (1+ \alpha) \bigr)$ and consequently $$ |N(A,b)|\le m^n
\bigl(6n^4 \log(1+ \alpha) \bigr) ^{n-1}.  \tag{cCook} $$ Thus, $|N(A,b)|\le
c_n m^n \log^{n-1}(1+ \alpha)$, where $c_n$ is some quantity depending only on
$n$.

Chirkov~\cite{ChirkovMA2} lowers the exponent of $m$ to $\antie{n/2}$
in the last estimate. Let
$$
v_0 , \dots , v_k, v_{k+1} , \dots , v_r
$$
be vectors in $\RR^n$ such that the system
$$
(v_0,1) , \dots , (v_k,1), (v_{k+1},0) , \dots , (v_r,0)
$$
is linear independent. The set
$$
\Conv\set{v_0 , \dots , v_k}+ \Cone\set{v_{k+1} , \dots , v_r}
$$
is called the {\em generalized simplex} with dimension $r$.
Denote by $\CharCone(P)$ the {\em cone of directions} of the polyhedron $P$
($\CharCone(P)$ consists of vectors $y$ such that for any $x \in
P$ and $\gamma > 0$ the point $x+ \gamma y$ belongs to $P$).
The set of generalized simlexes $S_1 , \dots ,
S_l$ is called the {\em covering} of $P$ if
$P=\cup _{i=1}^l S_i$ and for any $i=1 , \dots , l$
every vertex of $S_i$ and every edge of $\CharCone(S_i)$
are a vertex of $P$ and an edge of $\CharCone(P)$ respectively.  Let
$$
\xi _n (m) =
{m-\lfloor \frac{n-1}{2} \rfloor -1 \choose \lfloor \frac{n}{2} \rfloor
}+{m-\lfloor \frac{n}{2} \rfloor -1 \choose \lfloor \frac{n-1}{2} \rfloor }.
$$
As shown in~\cite{Charatheodory}, for a polyhedron $P$
there is a covering containing at most $n!\,\xi_n(M)$ generalized
simplexes.  Denote by $N_i$ the set of vertices of the convex hull of
integer points in $S_i$.  For the polytope $P$ we obviously have
$N(A,b) \subseteq \cup _{i=1}^l N_i$. Chirkov~\cite{ChirkovMA2}
bounded the quantity of coefficients in the linear system associated with $S_i$
in terms of $n$ and $\alpha$ and estimated $|N_i|$ using (\ref{cCook}).
These leaded to the following bound.
\begin{theorem} \label{Verh} \cite{ChirkovMA2}
$$ |N(A,b)|\le n^{7n}
\xi_n(m) \bigl(6 \log (1+ \alpha) + 3 \log n \bigr)^{n-1}.
$$
\end{theorem}
From Theorem~\ref{Verh} it follows that
$|N(A,b)|\le c_n m^{\antie{n/2}} \log^{n-1}(1+ \alpha)$,
where $c_n$ is some quantity depending only on $n$.  In Sect.\,3 we show
that for any fixed $n$ this estimate can not be improved with respect
to order.

\begin{remark}
It is known~\cite{McMullen} that the number of facets of a polyhedron $P$ in
$\RR^n$ with $v$ vertices is at most $\xi _n(v)$. From this and estimates above
it follows that for any fixed $n$ the number of facets in $P_I$ is bounded
above by some polynomial in $m$ and $\log (1+ \alpha)$.
This result and Lenstra's algorithm~\cite{Lenstra} allowed
Shevchenko~\cite{1984,book} to construct a {\em quasipolynomial} (i.e.
polynomial for fixed $n$) algorithm for finding all facets and vertices in
$P_I$.  \end{remark}

\begin{remark} Chirkov~\cite{ChirkovMA2} has obtained the upper bounds on the
number of vertices in the polytope $P_I$, as function of its {\em metric
characteristics}, i.e. the diameter $d= \max _{x,y \in P} \max _{j=1 , \dots ,
n} |x_j - y_j|$ and the volume $V$. It has be proved that for any fixed $n$ the
quantity $|N(A,b)|$ is bounded above by some polynomial in $m$ and $\log
(1+d)$, with its highest term being $c_n m^{\antie{n/2}} \log^{n-1}(1+ d)$.  If
$P_I$ has full affine dimension, then an analogous result holds also for the
volume of $P$.  \end{remark}

\section{Upper bounds independent of the right-hand sides}

In this section we consider methods for determination of upper bounds on
$|N(A,b)|$ independent of $b$. These bounds are useful if entries of a
vector $b$ are more greater than
$$
\alpha _1 = \max\set{|a_{ij}| (i=1 , \dots , m ; j=1 , \dots , n)}.
$$

Let the rank of $A$ be equal to $n$ and the solution set $P$ associated with
$Ax \le b$ be not empty. Then denote by $V(A, b)$ the set of vertices of $P$,
and by $\Delta(A)$ the maximum of the moduli of $n$th-order minors of $A$, and
by $a_i$ the $i$th row of $A$.

Following~\cite{book}, Sect.\,3.6.2, let
$$
J(v)= \set{i, b_i -a_iv <n\Delta(A)}, \quad
M(v) = \set{x \in \ZZ^n , a_ix \le b_i \mbox{ for } i \in J(v)}
$$
for every vertex $v$ of $P$. Denote by $N(v)$ the set of vertices of the convex
hull of $M(v)$. In~\cite{book} it has proved that
$$
N(A,b)\subseteq \bigcup_{v \in V(A,b)} N(v).
$$
Let $\lceil v \rfloor$ denote a vector obtained from $v$ by
rounding off the components of $v$.
As in~\cite{book}, we make the change of variables
$x=x'+\lceil v \rfloor$ and let $b_i '=b_i - a_i \lceil v \rfloor$.
The system of inequalities describing $M(v)$
passes into the system $a_i x' \le b_i '$, $i \in J(v)$, and,
since for $i \in J(v)$ it holds that
$$ |b_i'|=\bigl| (b_i-a_i v)+(a_i v - a_i
\lceil v \rfloor ) \bigr| < n\Delta (A) + \frac{n \alpha _1}{2},
$$
then Theorem~\ref{Verh} leads us to
\begin{theorem} \label{myimprovement}
$$
|N(A,b)| \le c_n m^{\antie{n/2}} \log ^{n-1} (1+ \alpha _1),
$$
where $c_n$ is some quantity depending only on $n$.
\end{theorem}

Let us now consider the knapsack polytope.  For $i=1 , \dots , n$ we
denote by $M_i$ the solution of solutions to the system
$$
\Sj a_j x_j \le a_0,\quad x_i \in \ZZ,\quad x_j \in \ZZ, \quad x_j \ge 0
\quad (j=1 , \dots , i-1,i+1 , \dots , n),
$$
and by $N_i$ the set of vertices of $\Conv M_i$, and set $N_0 = \set{0}$.
\begin{theorem} \label{knap} \cite{1981}
\begin{enumerate}
  \item $|N_i| \le \left( \antie {\displaystyle \log a_i}+n-1 \atop
        \displaystyle n-1 \right)$;
  \item if $a_0 \ge \alpha _1 (\alpha_1-1)$, then
        $\NNN(a,a_0)=\bigcup\limits_{i=0}^n N_i$ and
$$
|\NNN(a,a_0)| \le 1 + \Sj {\antie{\log a_i}+n-1 \choose n-1} \le 1+
n\bigl(\antie {\log \alpha_1} + 1 \bigr) ^{n-1}.
$$
\end{enumerate}
\end{theorem}

\section{Lower bounds on the number of vertices}

The natural approach for determination of lower bounds on the
number of vertices in integer linear programming problems is
constructing special examples.  Rubin's work~\cite{Rubin} was the first in
this area. Rubin found the class of polyhedra with arbitrarily many (but
finite) vertices and facets.  The $k$th polytope in the class has $k+3$
vertices; it is given by the inequality
$F_{2k}x+ F_{2k-1}y \le F_{2k+1}^2 -1$, where $F_s$ is the $s$th Fibonacci
number, i.e.  $F_1=F_2=1$, $F_s=F_{s-1}+F_{s-2}$ $(s=3,4,\dots )$.

Veselov and Shevchenko~\cite{26} also investigated the two-dimensional knapsack
problem. For naturals $a,b,c$, we denote by $\NNN(a,b,c)$ the set of vertices
of the knapsack polytope defined by the inequality $ax+by\le c$.
Let $\beta_1=1,
\beta_2=2, \gamma_2=\gamma'_2=1$, and, for every natural $s\ge 2$, let
$\beta_{s+1}=2\beta_s+\beta_{s-1}$, $\gamma_{s+1}=\beta_s+\gamma'_s$, $
\gamma'_{s+1}=\beta_{s+1}-\beta_s+\gamma_s$. We note that
$$
\beta_s = \frac{(1+\sqrt{2})^s - (1-\sqrt{2})^s}{2\sqrt{2}}.
$$

\begin{theorem} \label{2dknapsack} \cite{26}
If $a<\beta_{s-1}$, then $|\NNN(a,b,c)| <2s$.
If $a=\beta_{s-1}$, then there are $b$, $c$ such that
  for $b = \beta_s$ and
$$
c=\left\{
\begin{array}{ll}
\gamma _s (\beta _{s-1} +\beta _s) & \mbox{ for even } s, \\
\gamma _{s+1} \beta _{s-1} +\gamma _{s-1} \beta _s) & \mbox{ for odd } s
\end{array}
\right.
$$
it holds that $|\NNN(a,b,c)| =2s$.
\end{theorem}
Thus, $\beta_{s-1}$ is the minimal value of $a$ such that
there exist $b, c$ for which $\beta_{s-1}= a\le b \le c$ and
$|\NNN(a,b,c)| =2s$.
We remark that, by Theorem~\ref{2dknapsack}, $|\NNN(a,b,c)| \le 2
\log _{1+\sqrt{2}}(1+2\sqrt{2}a)$.
Let $b_s = \min \set{b, b\ge a, |\NNN(a,b,c)| \ge s}$.
Veselov and Shevchenko proved~\cite{26} that
$b_{2s}=\beta _s$, $b_{2s+1}=\beta _s + \beta _{s-1}$.
Thus, the upper bounds for the two-dimensional knapsack problem can not be
improved. Analogous results has been obtained~\cite{26} for square systems
of linear inequalities and for systems of congruences.

For $n=3$, Morgan~\cite{Morgan} constructed the class of polytopes with $m=5$
such that $|N(A,b)|$ grows proportionally to $\log ^2 \alpha$.
Let $P$ be a polytope, given by the system
$$
\begin{array}{c}
x\ge 0 , y\ge 0 , z\ge 0,\\
x+ \psi y + \theta z \le \nu,\\
x+ \theta y + \varphi z \le \nu,
\end{array}
$$
where $\theta=2\cos \left(\frac{2\pi}{7}\right)$,
$\varphi=2\cos \left(\frac{4\pi}{7}\right)$, $\psi=2\cos
\left(\frac{6\pi}{7}\right)$. It is shown~\cite{Morgan} that
$|N(A,b)| \ge \frac{1}{32} \log ^2 \nu$.

Another construction is offered in~\cite{Barany}. For any $n \ge 1$,
B\'ar\'any et al. constructed a polyhedron $P$ with $m=2n^2$ such that
$|N(A, b)| \ge c_n \varphi ^{n-1}$, where $c_n$ is some positive quantity
depending only on $n$, and $\varphi$ is the length of the expression for the
rational coefficients in a system describing $P$.  (From this it follows that
for any fixed $m$ and $n$ the bound~(\ref{Verh}) cannot be improved with
respect to order.  We note that a stronger result follows from presented
below earlier Veselov's work~\cite{VeselovLowerKnapsack}. Chirkov
established~\cite{ChirkovMA1} that the upper bound~(\ref{Verh}) cannot be
improved with respect to order even if only $n$ (but not $m$) is fixed.)

We present another example. As shown in~\cite{1981}, the vertex set
$\NNN(a,a_0)$ of the knapsack polytope with inequality
$$
2^{n-1}x_1 + \dots + 2 x_{n-1} + x_n \le 2^n -1
$$
consists of $2^n$ points. We note that formula~(\ref{knapsack1}) yields
$\NNN(a,a_0) \le n!$. This example shows that
upper bounds on $|\NNN(a,a_0)|$ and $|N(A,b)|$ must bear (have?) exponential
character in $n$.

Another approach, due to Veselov~\cite{VeselovLowerKnapsack},
consists in constructing a lower bound for the {\em mean} number of vertices
in an integer linear programming problem.
For each vector $a(a_0, a_1 , \dots , a_n)$ with integer non-negative
components such that $a_n =1$ and $0 \le a_i \le
\Delta -1$ $(i=0 , \dots , n-1)$, let us consider the set
$\MMM ' (a, \Delta)$ of solutions to the system
$$
\Sj a_j x_j \equiv a_0\ (\mod \Delta),
  \quad x_j \in \ZZ,
  \quad x_j \ge 0
  \quad (j = 1, \dots , n)
$$
and the set $\NNN ' (a, \Delta)$ of vertices of the convex hull of $\MMM ' (a,
\Delta)$.  Then $\varphi (\Delta) = \Delta ^{-n} \sum
|\NNN ' (a, \Delta)|$ is the mean number of vertices, where the summation is
over all $a$ in the range under consideration.

Veselov~\cite{VeselovLowerKnapsack} estimates the number
$\sigma (p)$ of $a$ such that $p \in \NNN ' (a,
\Delta)$; and finds a bound on $\Delta ^{-n} \Sum
\sigma (p)$, where the summation is over all $p=(p_1 , \dots , p_n) \in \ZZ^n$
such that $p_j \ge 0$ $(j=1 , \dots , n)$ and $(p_1+1)\cdot \dots
\cdot(p_n+1) \le \Delta$. It is obvious, that $\varphi(\Delta)=\Delta
^{-n} \Sum \sigma (p)$. The resulting bound has the following form:
$$
\varphi(\Delta) \ge \frac{c_n}{(n-1)^{n-1}}
\Bigl( \log \Delta - n - 2 - n \log(n-1) \Bigr)^{n-1}, \tag{phiDelta}
$$
where
$c_n = \Bigl(4n3^n\bigl( (n-1)!  \bigr)^2 \Bigr)^{-1}$.

Now let us consider the set $\varphi (\gamma)$ of vectors $(a_0 , \dots ,
a_{n-1}, \gamma)$ with integer non-negative components such that $a_i \le
\gamma$ $(i=1 , \dots , n-1)$, $\gamma (\gamma -1) \le a_0 <\gamma ^2$.
The knapsack polytope given by the inequality
$$
a_1 x_1 + \dots
a_{n-1}x_{n-1}+ \gamma x_n \le a_0
$$
is considered for any $(a_0, a) \in \varphi (\gamma)$, where $a=(a_1 , \dots ,
a_{n-1}, \gamma)$. The quantity $\psi (\gamma) = \gamma^{-n} \Sum |\NNN
(a, a_0)|$, where the summation is over all $(a_0,a)$ in $\varphi (\gamma)$,
is the mean number of vertices in such a polytope. Since
$a_0 \ge \gamma (\gamma -1)$ and $a_i \le \gamma$ $(i = 1 , \dots , n-1)$,
then, by Theorem~\ref{knap}, $|\NNN(a, a_0)| \ge |N'|$, where $N'$ is the set
of vertices of the convex hull of integer solutions to the system
$$
\begin{array}{c}
a_1 x_1 + \dots + a_{n-1} x_{n-1} + \gamma x_n \le a_0,\\
x_j \ge 0 \quad (j=1 , \dots , n-1).
\end{array}
$$
It is obvious that $N'$ is mapped one-to-one into the set of vertices
of the convex hull of non-negative solutions to the congruence
$$
x_0 + a_1 x_1 + \dots a_{n-1} x_{n-1} \equiv a_0 \,(\mod \gamma),
$$
or, equivalently,
$$
x_0 + a_1 x_1 + \dots a_{n-1} x_{n-1} \equiv a_0 - \gamma(\gamma -1)\,(\mod
\gamma). $$
Thus, $\psi (\gamma) = \varphi (\gamma)$. In particular, this gives us
\begin{theorem}
\cite{VeselovLowerKnapsack}
There exist non-negative $a=(a_1 , \dots , a_n) \in \ZZ^n$ and $a_0 \in
\ZZ$ such that $a_j \le \alpha$ $(j=0 , \dots , n)$ and
$$
|\NNN (a_0,a)| \ge c_n \log ^{n-1} \alpha,
$$
where $c_n$ is some positive quantity depending only on $n$.
\end{theorem}
This result shows that for any fixed $n$ the order in (\ref{knapsack2}) can not
be decrease.

The treated approach was extended by Chirkov and
Shevchenko~\cite{ChirkovLower,ChirkovMA1,ChirkovKorsh} to
the case of an arbitrary system
$Ax \le b$. The following theorem holds.
\begin{theorem}
\cite{ChirkovMA1,ChirkovKorsh}
For any $n$, $m$ and $\alpha$ there exist a polyhedron such that
$$
|N(A,b)|\ge \frac{\xi_n(m)}{4^{n+2}n^n(n-1)!}  \ln ^{n-1} \alpha.
$$
\end{theorem}
Thus, the lower bound on the number of vertices has the following form:
$$
|N(A,b)|\ge c_n m^{\antie{n/2}} \log ^{n-1} \alpha,
$$
where $c_n$ is some positive quantity depending only on $n$.  This result shows
that for any fixed $n$ the bound (\ref{Verh}) cannot be improved.

\section{Final notes and conclusions}

Let us sum up.
For the number of vertices in the convex hull of integer points in a polyhedron
the following bound has been obtained~\cite{ChirkovMA2}:
$$
|N(A,b)| \le c_n m^{\antie{n/2}} \log ^{n-1} (1+\alpha),
\tag{UB}
$$
where $c_n$ is some quantity depending only on $n$. We have shown that
this result is valid even
if
$$\alpha = \max \set{|a_{ij}| (i=1 , \dots , m; j=1, \dots , n)}.$$
It is established~\cite{ChirkovMA1} that
this bound can not be improved for any fixed $n$ with respect to order, namely,
for any $n$, $m$ and $\alpha$ there exist a polyhedron such that
$$
|N(A,b)|\ge c'_n m^{\antie{n/2}} \log ^{n-1} \alpha,
$$
where $c'_n$ is some positive quantity depending only on $n$.  Analogous result
is also valid for the knapsack problem~\cite{VeselovLowerKnapsack}.

The proofs of these lower bounds are not constructive. Examples of polyhedra
with the ``large'' number of vertices in $P_I$ are known only for few values of
$m$. In particular, for $n \to \infty$ and for any positive $c_n$ no
sequence of knapsack polyhedrons with at least $c'_n \log^{n-1} \alpha$
vertices is known.

From the upper bound (\ref{UB}) on the number of vertices in $P_I$ and from the
inequality~\cite{McMullen}, estimating the number of facets in a
polyhedron with $v$ vertices, follows that the number of facets in $P_I$ is
at most
$$
c_n m^{{\antie{n/2}}^2} \log^{(n-1)\antie{n/2}}(1+\alpha).
$$
The problem of an attainability of this bound is open.

\bibliographystyle{alpha}
\bibliography{extrengl}
\addcontentsline{toc}{section}{References}

\end{document}